\documentclass[a4paper,11pt,twoside]{article}
\usepackage{fancyhdr,latexsym,amsmath,amsfonts,amssymb,amsbsy,amsthm,url}
\usepackage[hscale=0.7,centering]{geometry}
\usepackage{graphicx,epsfig}
\usepackage{amssymb,amsmath}
\usepackage{times}
\usepackage{color}
\usepackage{amsfonts,dsfont}[11pt]
\usepackage{amsmath}[11pt]
\newcommand{\X}{\mathds{X}}

\numberwithin{equation}{section}

\def\qed{ \hfill $\Box$}

\fancypagestyle{usual}{
  \lhead[\thepage]{}
  \chead[{\it \shortauthors}]{\sc \shorttitle}
  \rhead[]{\thepage}
  \lfoot[]{}
  \cfoot[]{}
  \rfoot[]{}
  
}

\makeatletter

\renewcommand{\section}{
  \@startsection
  {section}
  {1}
  {0pt}
  {1.1\baselineskip}
  {0.2\baselineskip}
  {\sc \centering}
}
\makeatother

\newcommand{\Z}{\mathbb Z}
\newcommand{\veps}{\varepsilon}

\definecolor{db}{rgb}{0.1,0,0.75}
\definecolor{lm}{cmyk}{0 ,1,0,0}

\newcommand{\R}{\mathbb R}

\newcommand{\pr}{\mathbb P}
\newcommand{\N}{\mathbb N}

\newtheorem{lem}{Lemma}

\newtheorem{prop}{Proposition}
\newtheorem{thm}{Theorem}

  \DeclareMathOperator{\E}{E}

\DeclareMathOperator{\Tr}{Tr}

\newcommand{\I}{\mathbb{I}}

\theoremstyle{definition} \theoremstyle{remark}

\begin{document}
\def\shorttitle{Balanced matrices}
\def\shortauthors{A. Basak,  A. Bose, A. Sen}
\title{\sc Balanced random Toeplitz and Hankel Matrices}
\normalsize
\author{
\parbox[t]{0.45\textwidth}{{\sc Anirban Basak}\thanks{Anirban Basak is supported by Melvin and Joan Lane Endowed Stanford Graduate Fellowship Fund.}\\
Department of Statistics\\
Stanford University\\
USA
}
\parbox[t]{0.35\textwidth}{{\sc Arup Bose}\thanks{ Research
supported by J.C.Bose Fellowship, Department of Science and Technology, Government of India.}\\
Statistics and Mathematics Unit\\
Indian Statistical Institute\\
202 B. T. Road,
Kolkata 700108\\
INDIA
 }
 }
\date{January 10,  2010} \maketitle
\begin{abstract} Except the Toeplitz and Hankel matrices, the common patterned matrices for which the limiting spectral distribution (LSD) are known to exist,  share a common property--the number of times each random variable appears in the matrix is (more or less) same across the variables.
 Thus it seems natural to ask what happens to the spectrum of the Toeplitz and Hankel matrices when each entry is scaled by the square root of the number of times that entry appears in the matrix instead of the uniform scaling by
$n^{-1/2}$.  We show that the LSD of these balanced matrices exist and derive integral formulae for the moments of the limit distribution. Curiously, it is not clear if these moments define a unique distribution.
\end{abstract}


\section{Introduction and main results}\label{section: top_han_var}
For any (random and symmetric)   $n\times n$ matrix $B$, let $\mu_1(B), \ldots , \mu_n(B) \in
\mathbb R$ denote its eigenvalues including multiplicities. Then the
empirical spectral distribution (ESD) of $B$ is the (random) distribution function on
$\mathbb R$ given by
$$F^{B} (x) = n^{-1}  \# \Big \{ j :  \mu_j(B)  \in (-\infty, x], \ 1 \le j \le n  \Big \}.$$  \noindent For a sequence of random $n \times n$ matrices $\{ B_n\}_{ n \ge 1}$ if as $n \to \infty$,
the corresponding ESDs $F^{B_n}$ converge  weakly (either almost surely or in probability) to a (nonrandom) distribution
$F$ in the space of probability measures on $\mathbb R$, then $F$ is called the limiting spectral distribution (LSD) of $\{B_n\}_{n \ge 1}$.  See Bai (1999)\cite{Bai99}, Bose and Sen (2007)\cite{bosesen08} and Bose, Sen and Gangopadhyay (2009)\cite{bosesengangopadhyay09} for description of several interesting patterned matrices whose  LSD exists. Examples include the Wigner, the circulants, the Hankel and the Toeplitz matrices.
For  the Wigner and circulant matrices, the number of times each random variable appears in the matrix is (more or less) same across the variables. We may call them balanced matrices. For them, the LSD exist after the eigenvalues are scaled by $n^{-1/2}$.

\begin{figure}[htbp]
 \begin{center}
 \epsfig{file=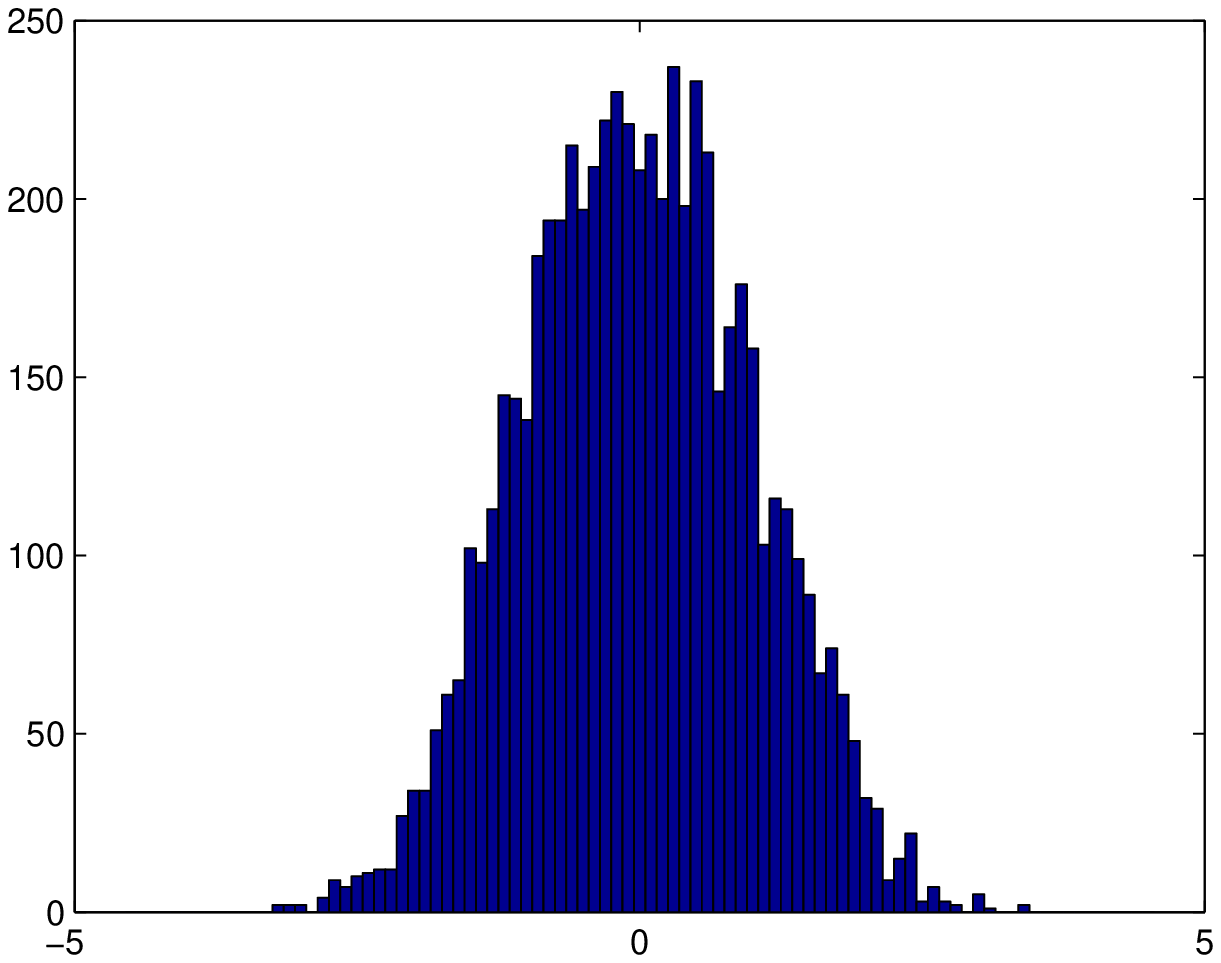,width=2.25in,height=1.1in}
 \epsfig{file=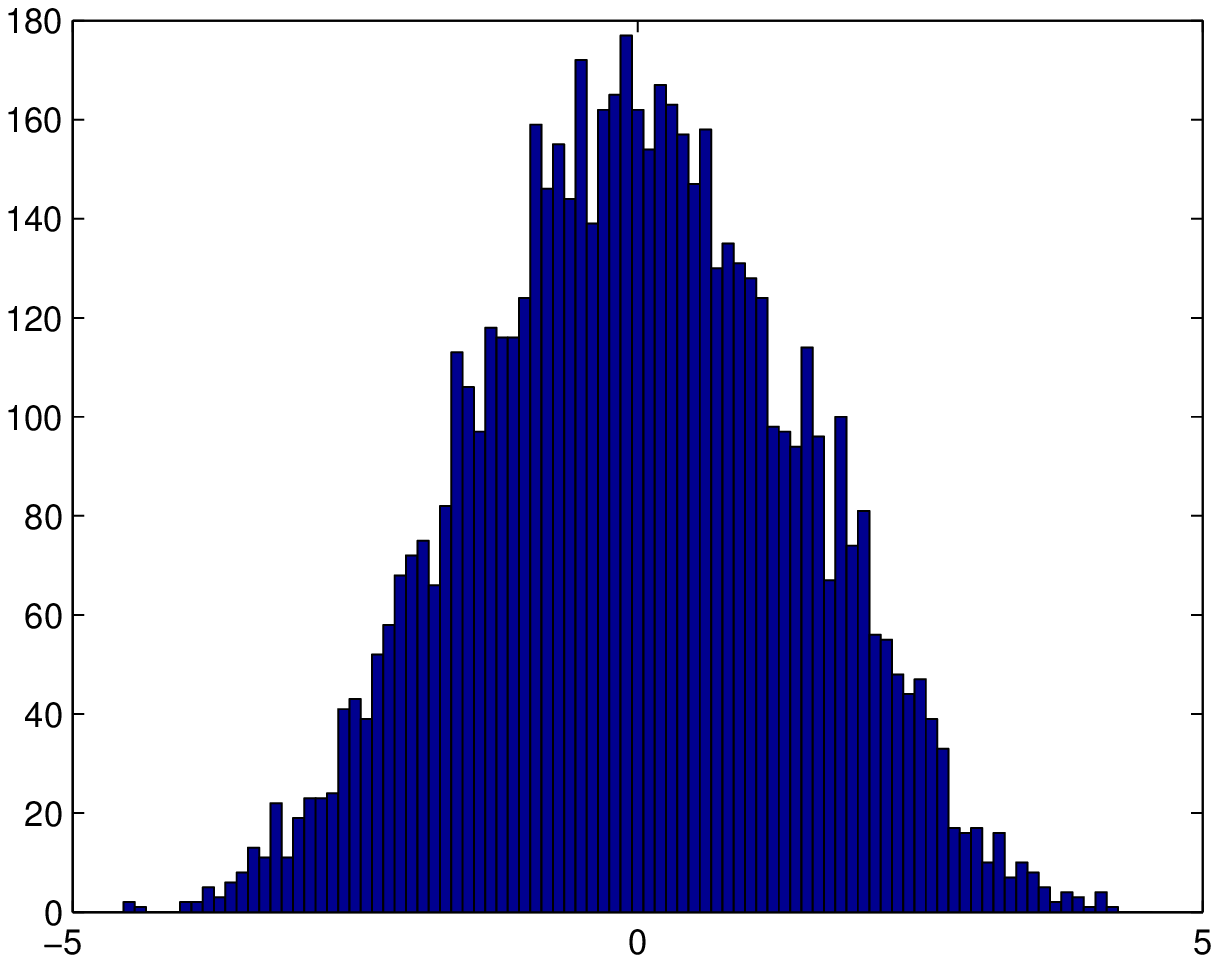,width=2.25in,height=1.1in}
\epsfig{file=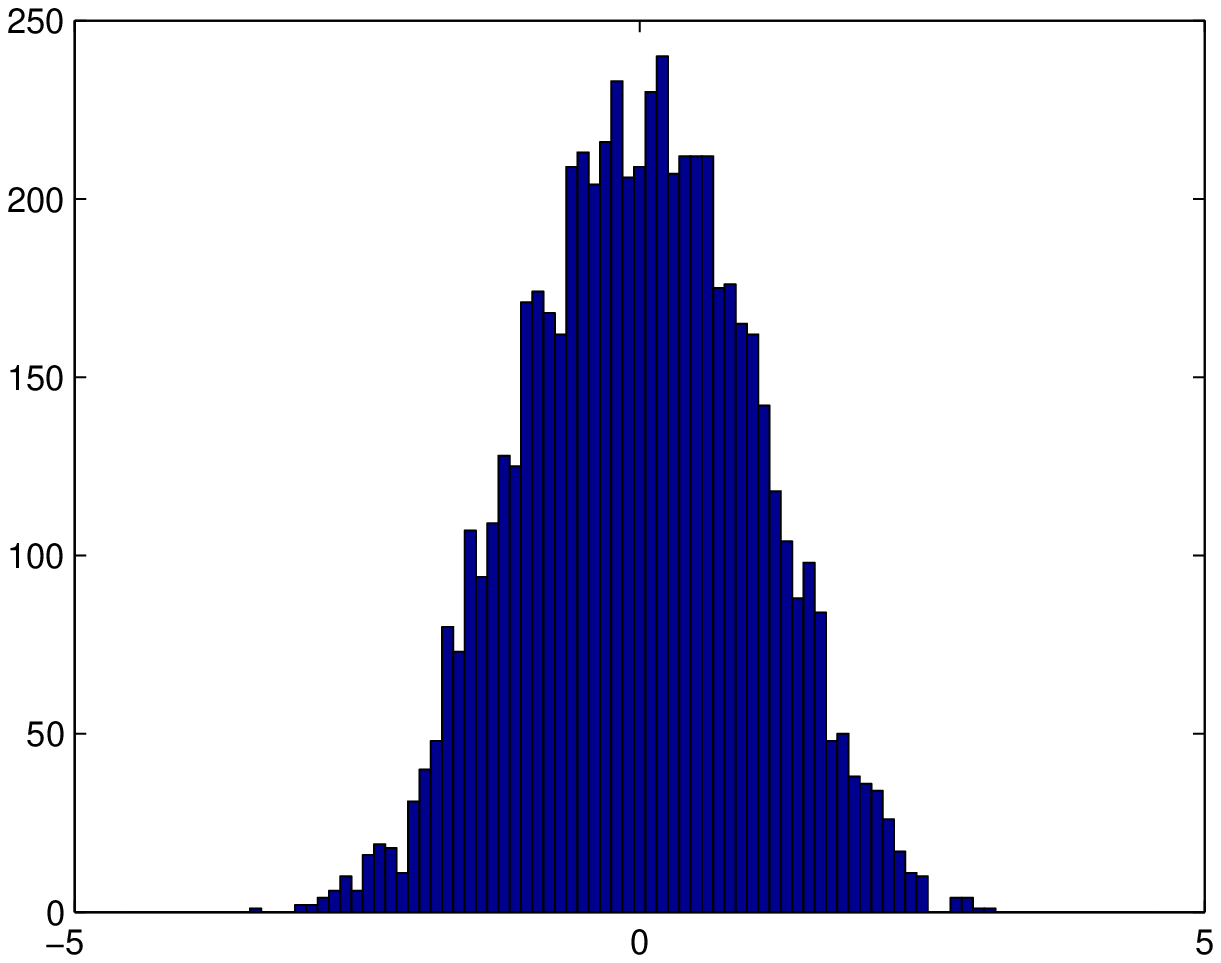,width=2.25 in,height=1.1in}
\epsfig{file=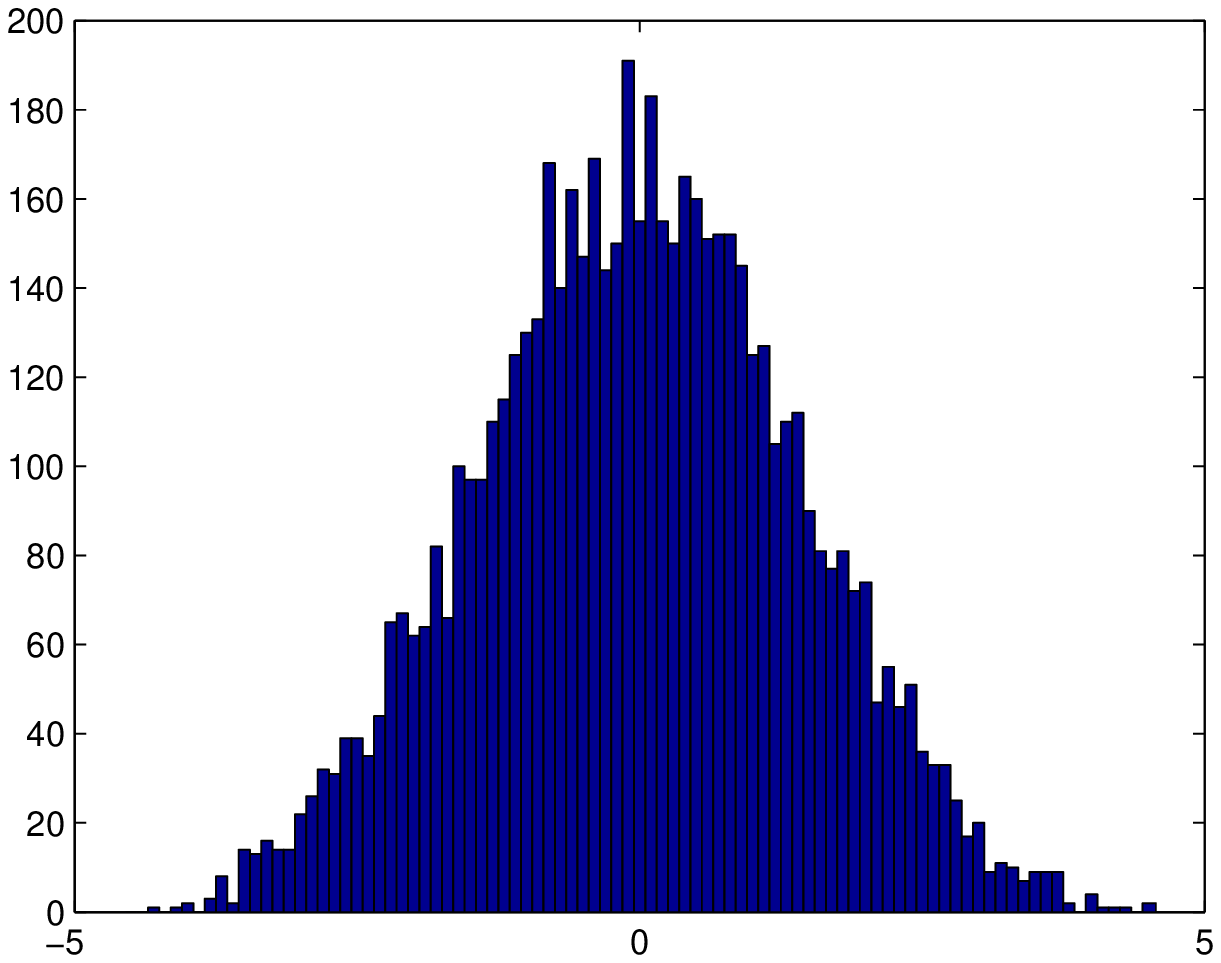,width=2.25 in,height=1.1in}
 \caption{\small Histograms of the ESD of $15$ realizations of the Toeplitz matrix (left) and the balanced Toeplitz matrix (right)  of order $400$ with standardized Normal$(0,1)$ (\texttt{top row}), 
 and Bernoulli$(0.5)$ (\texttt{bottom row})
entries.}
\label{fig:figtoep}
 \end{center}
 \end{figure}

\begin{figure}[htbp]
 \begin{center}
 \epsfig{file=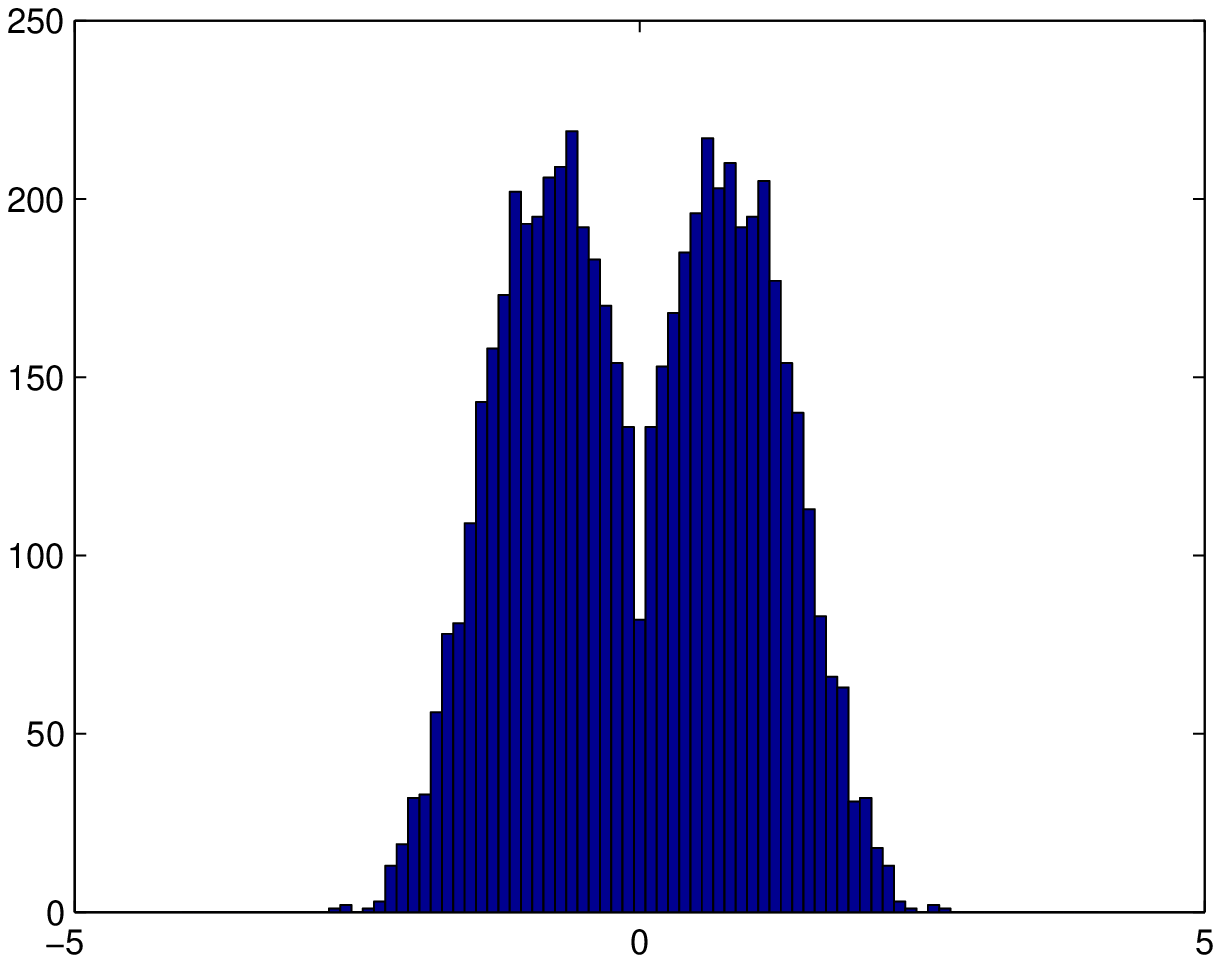,width=2.5in,height=1.1in}
 \epsfig{file=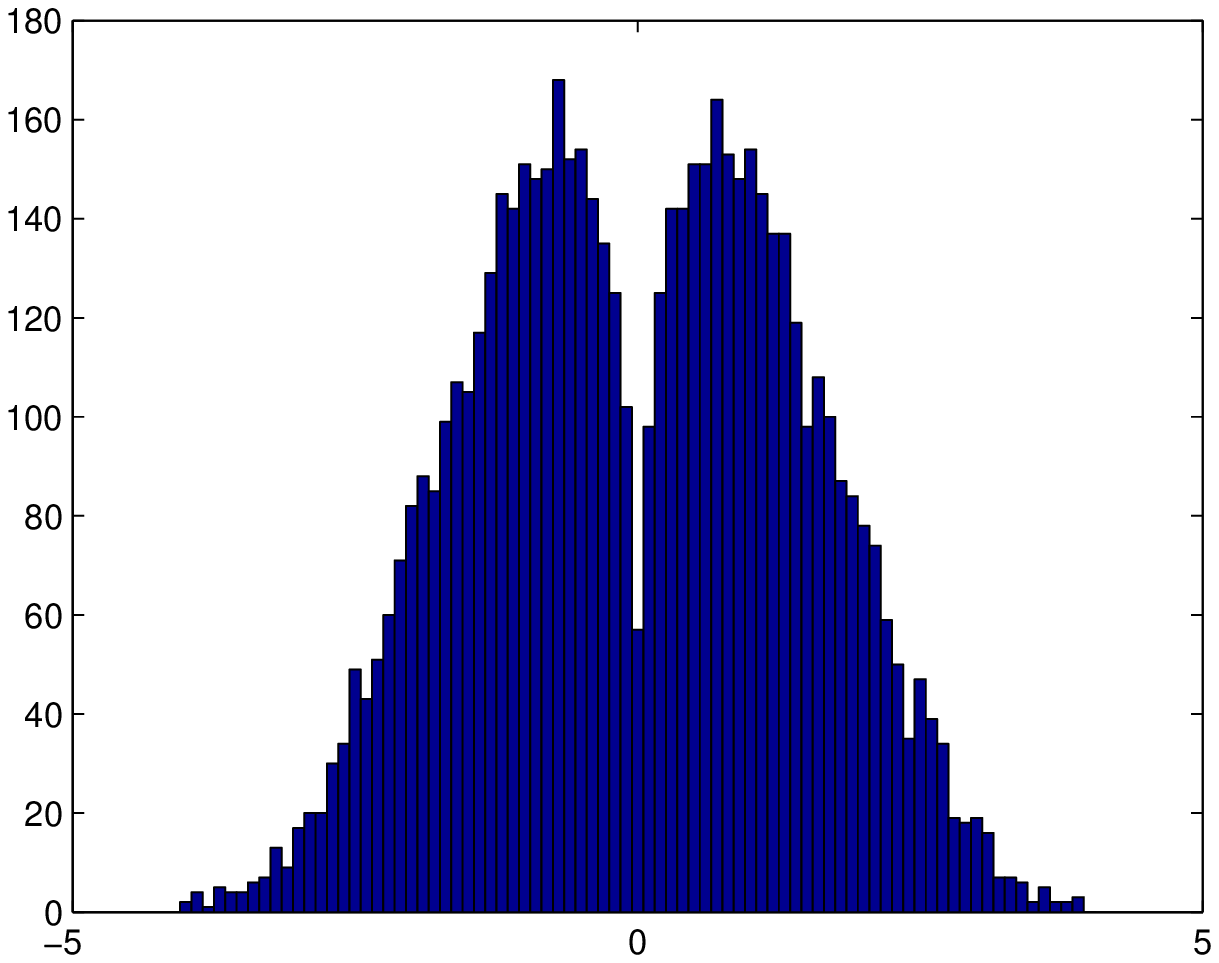,width=2.5in,height=1.1in}
\epsfig{file=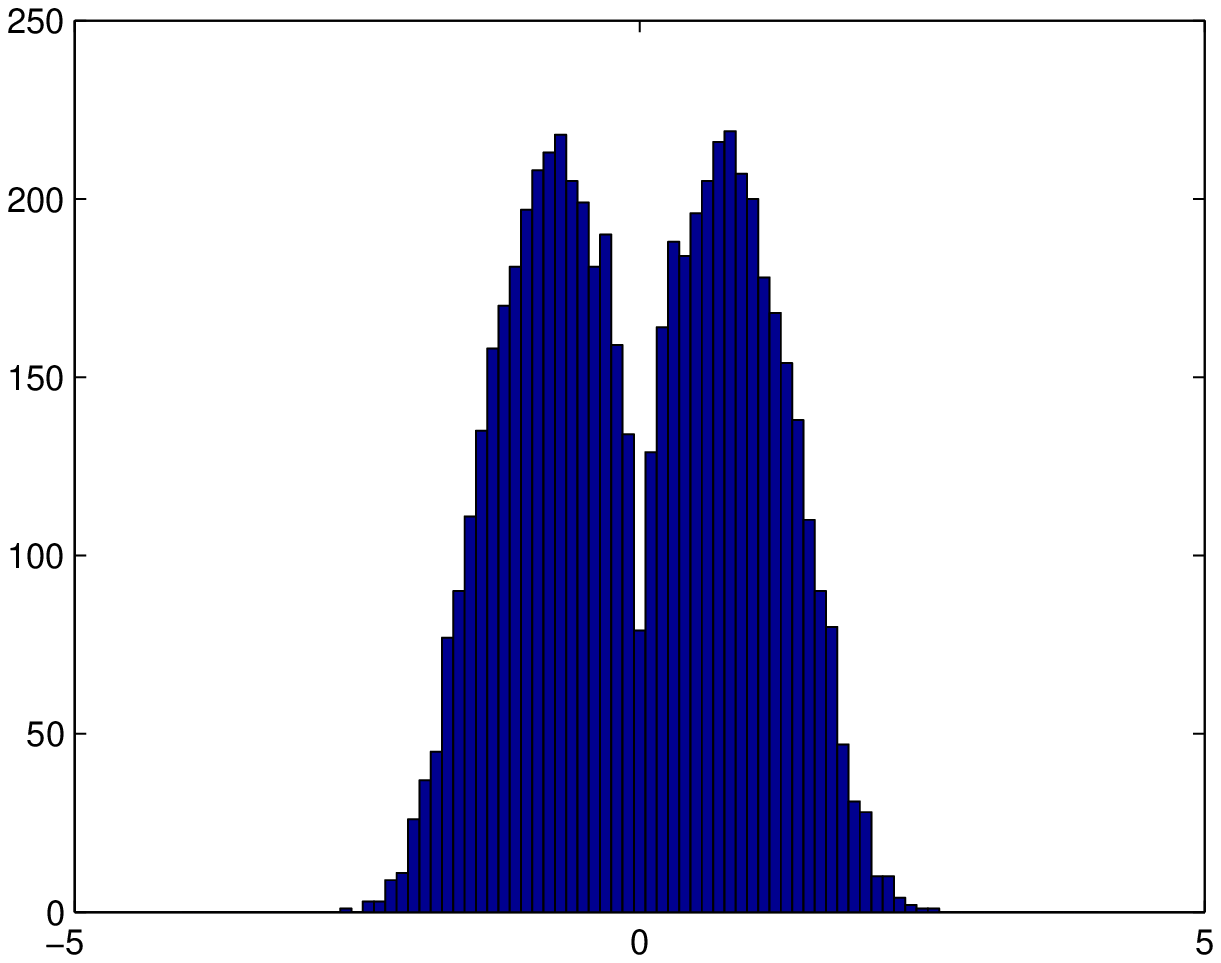,width=2.5 in,height=1.1in}
\epsfig{file=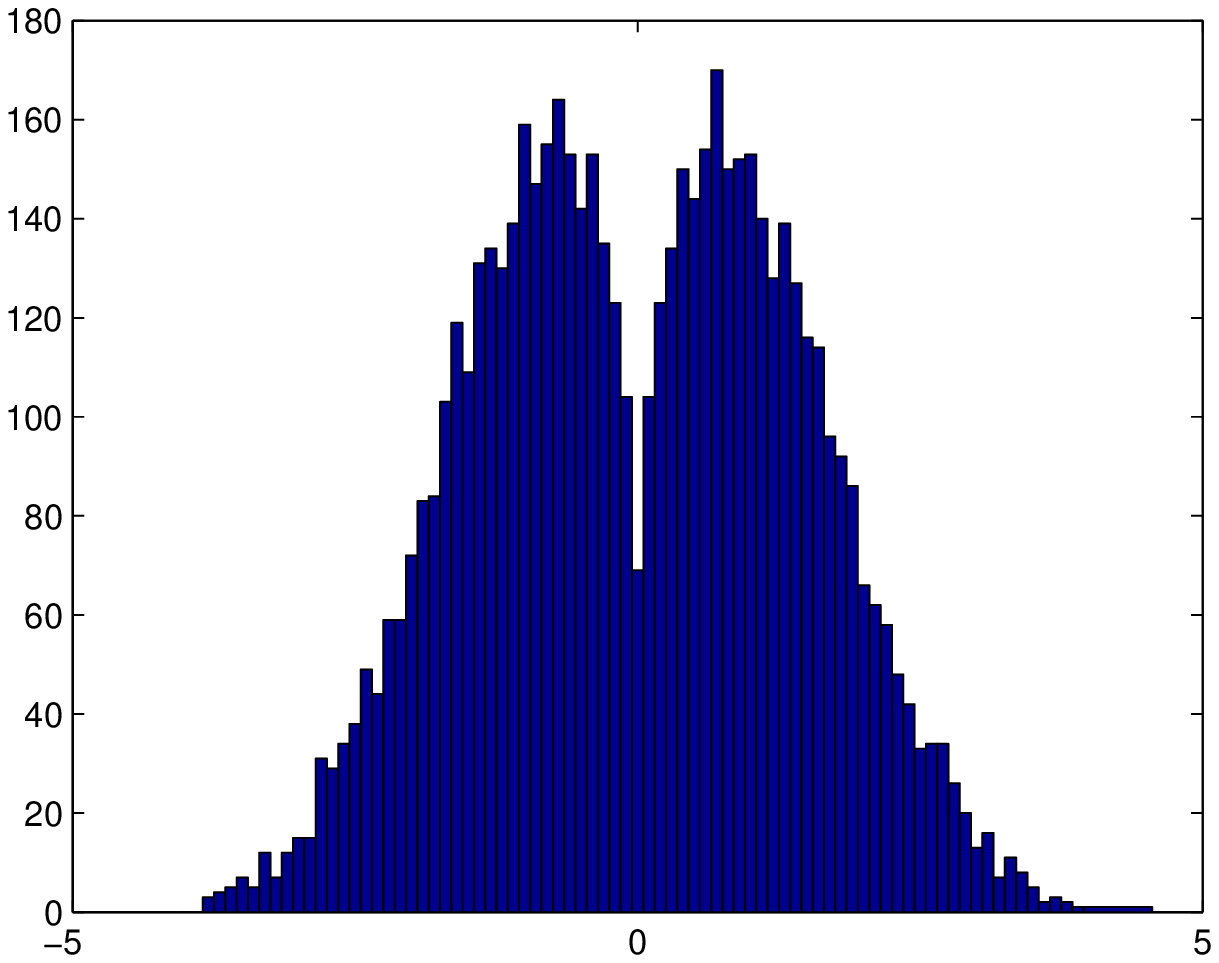,width=2.5 in,height=1.1in}
\caption{\small Histograms of the ESD of  $15$ realizations of the Hankel matrix (left) and the balanced Hankel matrix (right)  of order $400$ with standardized Normal$(0,1)$ (\texttt{top row}), 
 and Bernoulli$(0.5)$ (\texttt{bottom row})
entries.}
\label{fig:fighankel}
 \end{center}
 \end{figure}

 Consider the $n \times n$ symmetric Toeplitz and Hankel matrices with an i.i.d. \textit{input} sequence $\{x_i\}$. For these matrices when scaled by $n^{-1/2}$, the LSD  exists.
 The limit is non Gaussian   but not much more is known about its properties. See Bryc, Dembo and Jiang (2006)\cite{bryc} and Hammond and Miller (2005)\cite{hammil05}. However, these matrices are  unbalanced.
It seems natural to consider a balanced version of the Toeplitz and Hankel matrices where each entry is scaled by the square root of the number of times that entry appears in the matrix instead of the uniform scaling by $n^{-1/2}$.
Define the (symmetric) balanced  Hankel  and Toeplitz matrices
$BH_{n}$ and $BT_{n}$ with input $\{x_i\}$  as follows:
\begin{equation}\label{def:hankel}
BH_{n}\  = \left[ \begin{array} {cccccc}
            \frac{x_{1}}{\sqrt{1}} & \frac{x_{2}}{\sqrt{2}} & \frac{x_{3}}{\sqrt{3}} & \ldots & \frac{x_{n-1}}{\sqrt{n-1}} & \frac{x_{n}}{\sqrt{n}} \\
            \frac{x_{2}}{\sqrt{2}} & \frac{x_{3}}{\sqrt{3}} & \frac{x_{4}}{\sqrt{4}} & \ldots & \frac{x_{n}}{\sqrt{n}} & \frac{x_{n+1}}{\sqrt{n-1}} \\
            \frac{x_{3}}{\sqrt{3}} & \frac{x_{4}}{\sqrt{4}} & \frac{x_{5}}{\sqrt{5}} & \ldots & \frac{x_{n+1}}{\sqrt{n-1}} & \frac{x_{n+2}}{\sqrt{n-2}} \\
                  &       &       & \vdots &       &   \\
            \frac{x_{n}}{\sqrt{n}} & \frac{x_{n+1}}{\sqrt{n-1}} & \frac{x_{n+2}}{\sqrt{n-2}} & \ldots & \frac{x_{2n-2}}{\sqrt{2}} & \frac{x_{2n-1}}{\sqrt{1}}
            \end{array} \right].
\end{equation}
\begin{equation}\label{def:toep_var}
    BT_{n}\  = \left[ \begin{array} {cccccc}
            \frac{x_{0}}{\sqrt{n}} & \frac{x_{1}}{\sqrt{n -1}} & \frac{x_{2}}{\sqrt{n - 2}} & \ldots & \frac{x_{n-2}}{\sqrt{2}} & \frac{x_{n-1}}{\sqrt{1}} \\
            \frac{x_{1}}{\sqrt{n -1}} & \frac{x_{0}}{\sqrt{n}} & \frac{x_{1}}{\sqrt{n-1}} & \ldots & \frac{x_{n-3}}{\sqrt{3}} & \frac{x_{n-2}}{\sqrt{2}} \\
            \frac{x_{2}}{\sqrt{n-2}} & \frac{x_{1}}{\sqrt{n-1}} & \frac{x_{0}}{\sqrt{n}} & \ldots & \frac{x_{n-4}}{\sqrt{4}} & \frac{x_{n-3}}{\sqrt{3}} \\
                  &       &       & \vdots &       &   \\
            \frac{x_{n-1}}{\sqrt{1}} & \frac{x_{n-2}}{\sqrt{2}} & \frac{x_{n-3}}{\sqrt{3}} & \ldots & \frac{x_{1}}{\sqrt{n-1}} & \frac{x_{0}}{\sqrt{n}}
            \end{array} \right].
\end{equation}
Strictly speaking $BT_n$ is not completely balanced, with the main diagonal being unbalanced compared to the rest of the matrix but this does not affect the asymptotic behavior of the eigenvalues. We use the above version because of the convenience in writing out the calculations later.
Figures \ref{fig:figtoep} and \ref{fig:fighankel} exhibit the  simulation results for the ESD of the above matrices.
We prove the following theorem.
\begin{thm}
\label{thm:main theorem} Suppose $\{x_i\}$ are independent mean zero variance one random variables which are either  uniformly  bounded, or are i.i.d. Then almost surely  the LSD, say $BT$ and $BH$  of the matrices $BT_{n}$ and $BH_{n}$ respectively, exist.
\end{thm}
\noindent\textbf{Remark} The limits distributions have unbounded support, are symmetric about zero and have all moments finite. The integral formulae for these moments  are given in \ref{eq:limexptoep} and \ref{eq:limexphan} later in Section \ref{subsection:conv} after we develop the requisite notation to write them out. Both LSD are non Gaussian and do not depend on the underlying distribution of the input sequence. It does not seem to be apparent from our formulae if these moments define a distribution uniquely. Establishing further properties of the limits is a difficult problem.
\section{Proof of Theorem \ref{thm:main theorem}}
\noindent
(1) In Section \ref{subsection: unif_bound} we first show that we may restrict attention to bounded $\{x_i\}$. 
  
\noindent 
 (2) In Sections \ref{subsection:momentsandtrace} --\ref{subsection: reduction_2} we develop the trace formula for moments, some related notions and results to reduce the  number of terms in the trace formula. In Section \ref{subsection:conv} we show that the \textit{expected} moments of the ESD of $\{BT_n\}$ converge. However, it does not seem to be straightforward to show that this limiting sequence uniquely determines a distribution. Even if it did, it is not clear how the convergence of the expected moments can be sharpened to convergence of the ESD itself. If we pull out the usual scaling $n^{-1/2}$, the  scaling for the $(i,j)$th entry is  $[1-\frac{|i-j|}{n}]^{-1/2}$ whose maximum is $n^{1/2}$. This unboundedness creates problems in the usual argument.
 
 \noindent 
  (3) In Section \ref{subsection:approximation} we discuss a known approximation result. 
  
  \noindent (4) Fix any $\varepsilon >0$. Let $BT_{n}^{\varepsilon}$ denote the top-left $n(1-\varepsilon) \times n(1- \varepsilon)$ principal sub-matrix of $BT_{n}$. The L\'{e}vy distance between $F^{BT_{n}}$ and $F^{BT_{n}^{\varepsilon}}$ is less than $\varepsilon$. Since these truncated matrices are well behaved, convergence of $F^{BT_{n}^{\varepsilon}}$ to a non-random distribution $BT^{\varepsilon}$ almost surely,  and also the corresponding expected moments follow by the same arguments as for the usual Toeplitz matrices.
This limit is uniquely determined by its moments.  This is done in Section \ref{subsection:limitoftruncation}.
  
\noindent (5) In Section \ref{subsection:connect}, we show by using results derived in (3) that as $\varepsilon \to 0$, the spectral measures of $BT^{\varepsilon}$ converge to some $F^{BT}$ and this is the LSD of $\{BT_{n}\}$. We finally use a uniform integrability argument to conclude that the moments of $F^{BT}$ are same as those obtained in  (2) above.

A similar proof works for the balanced Hankel by defining the truncated Hankel matrix obtained deleting the first $n \veps/2$ and last $n \veps/2$ rows and columns from $BH_n$.
\subsection{Reduction to uniform bounded case}
\label{subsection: unif_bound}
\begin{lem}
\label{lem:truncation} Suppose for every bounded, mean zero and variance one i.i.d. input
sequence $\{ x_{0} , x_{1} , x_{2} , \ldots \}$, $\{F^{BT_{n}}\}$ converges to some
non-random
distribution $F$ a.s. Then the same limit continues to hold if $\{x_i\}$
is i.i.d. with
mean zero and variance one. All the above hold for $\{F^{BH_{n}}\}$.
\end{lem}
We  make use of the bounded Lipschitz metric. It is defined on the space of probability measures as:
$$d_{BL}(\mu, \ \nu) = \sup \{ \int f d\mu - \int f d\nu : ||f||_{\infty} + ||f||_L \le 1\} $$
where $||f||_{\infty} = \sup_x |f(x)|, \ ||f||_L = \sup_{x \ne y}
|f(x)-f(y)|/|x-y|$.
Recall that convergence in $d_{BL}$ implies the weak convergence of measures and vice versa.

We also need the following fact. This Fact is an estimate of the metric
distance $d_{BL}$ in terms of trace. A  proof may be found in Bai and Silverstein (2006)\cite{baisilverstein06} or Bai (1999)\cite{Bai99}.

\noindent
{\bf Fact 1.}   Suppose $A, B$ are $n \times n$ symmetric real
matrices. Then
\begin{equation}\label{blmetric} d^2_{BL}(F^A, F^B) \le \left (\frac{1}{n}
\sum_{i=1}^n | \lambda_i(A) - \lambda_i(B)| \right)^2 \le \frac{1}{n}
\sum_{i=1}^n ( \lambda_i(A) - \lambda_i(B)) ^2 \le \frac{1}{n} \Tr(A
-B)^2.
\end{equation}

\noindent
\textbf{Proof of Lemma \ref{lem:truncation}}. For brevity, we deal with
only the balanced Toeplitz.
Same arguments work for the balanced Hankel matrix.
\noindent
For $t > 0,$ define,
$$\mu (t) \stackrel{def}{=} E [ x_{0} (I | x_{0} | \leq t ) ], \ \
\sigma^{2} (t) \stackrel{def} {=} Var ( x_{0} I( | x_{0} | \leq t) ) = E[ x_{0}^ {2} I | x_{0}|
\leq t ) ] - \mu ( t) ^ {2},$$
$$x_{i}^{*} = \frac{x_{i} I ( | x_{i} | \leq t ) - \mu(t) } { \sigma(t) } = \frac{ x_{i} - \bar{ x}_{i} } { \sigma (t) },\   \text{where}\  \bar{x}_{i} = x_{i} I ( | x_{i} | > t ) + \mu ( t) = x_{i} - \sigma ( t) x_{i} ^ {*}.$$

 Let $\{ BT_{n}^{*} \}$ be the balanced Toeplitz matrix for the
input sequence $\{ x_{i}^{*} \} $ and $\{\widetilde{\overline{ BT}}_{n} \}$ be the same for the input sequence $
\{\bar{x}_{i} \} $. It is clear that $\{ x_{i}^{*} \} $ is a bounded, mean zero, variance one
i.i.d. sequence. Hence by our assumption, $F^{BT_{n}^{*}}$ converges to a non-random distribution
function $F$ a.s. Using Fact 1,
\begin{eqnarray*}
d^{2}_{BL} ( F^ {BT_{n}}  , F^{BT_{n}^{*}} ) & \leq & 2 d^{2}_{BL} ( F^{BT_{n}}, F^{ \sigma (t)
BT_{n}^{*} } ) + 2 d^{2}_{BL} ( F^{BT_{n}^{*}} , F^{\sigma (t) BT_{n}^{*} } )  \\
& \leq &\frac{2}{n}
Tr [ ( BT_{n}- \sigma (t) BT_{n}^{*} ) ^{2} ] + \frac{2}{n} (1- \sigma (t) ) ^{2} Tr [(BT_{n}^{*})^{2}
].
\end{eqnarray*}
Now using strong law of large numbers, we get,
\begin{eqnarray*}
\frac{1}{n} Tr[ (BT_{n}^{*} )^{2} ]
& = &\frac{1}{n}\sum_ {i,j} \left( \frac{x_{|i -j |}^{*}} { \sqrt{n - |i -j| } } \right) ^{2}\\
& = & \frac{1}{n}(n \times \frac{{x_{0}^{*}}^{2}}{n}   +2(n-1) \times \frac{{x_{1}^{*}}^{2}}{(n-1)}
+ \cdots +  2 \times \frac{{x_{n-1}^{*}}^{2}}{1} )\\ &\leq &\frac{2}{n}({x_{0}^{*}}^{2} +
{x_{1}^{*}}^{2} +\cdots + {x_{n-1}^{*}}^{2} ) \stackrel{a.s.}{\rightarrow} 2E({x_{0}^{*}}^{2}) = 2.
\end{eqnarray*}
Note that, \ $1 - \sigma(t) \rightarrow 0 $ as $ t \rightarrow \infty.$
Similarly,
\begin{eqnarray*}
\frac{1}{n} Tr[ ( BT_{n}- \sigma (t) BT_{n}^{*} ) ^{2} ] & = &\frac{1}{n} Tr[ \widetilde{\overline{BT}}_{n}^{2}]\\
& = &
\frac{1}{n} \sum_{i,j} \left( \frac{\bar{x} _{ |i-j|}
} { \sqrt{n - |i-j|} } \right) ^{2} \\
& = & \frac{1}{n} ( n \times \frac{\bar{x}_{0}^{2}} {n} + 2(n-1) \times \frac{ \bar{x}_{1}^{2} }
{(n-1) } + \cdots + 2 \times \frac{\bar{x}_{n-1}^{2} } {1} )\\
& \leq  & \frac{2}{n}( \bar{x} _{0}^{2} + \bar{x}_{1}^{2} + \cdots + \bar{x}_{n-1}^{2} )
 \stackrel {a.s.} {\rightarrow} 2 E[ \bar{x} _{0}^{2} ]=
 1 - 2\mu ( t) ^{2} -\sigma^{2} (t)  \rightarrow 0\ \  \  \text{as} \ \   t \rightarrow \infty.
 \end{eqnarray*}
Hence combining the above arguments, we get,
$\limsup_{n} d_{BL}( F^{BT_{n}} , F^{BT_{n}^{*} } ) \rightarrow 0 \ \text{a.s. as} \ \
t \rightarrow \infty.$
This completes the proof of this Lemma. \qed

\subsection{Moment and trace formula}\label{subsection:momentsandtrace}
We need some notation to express the moments of the ESD in a way convenient to further analysis.

\noindent
{\bf Circuit and vertices:} A \emph{circuit} is any function $\pi: \{0,1,2,\ldots,h\} \rightarrow \{1,2,\ldots,n\}$ such that $\pi(0)=\pi(h)$. Any $\pi(i)$ is a \emph{vertex}.
A circuit depends on $h$ and $n$ but we will suppress this dependence.

Define two functions $L^T$ and $L^H$, which we call \textit{link functions}, and as
\begin{equation}
L^{T}(i,j)=|i-j|\ \ \text{and} \ \ L^{H}(i,j)=i+j-1.
\end{equation}
Also for $L=L^H$ or $L^T$, as the case may be, define
$$\X_{\pi} = x_{L(\pi(0), \pi(1))} x_{L(\pi(1),\pi(2))} \cdots
x_{L(\pi(h-2), \pi(h-1))} x_{L(\pi(h-1), \pi(h))}.$$
Also define the following
\begin{equation}
\phi_{T}(i,j)=n - |i-j|\ \ \text{and} \ \ \phi_{H}(i,j)=\min (i+j-1,2n-i-j+1).
\end{equation}
\begin{equation}
\phi_T^n(x,y)=\phi_T^{\infty}(x,y)= 1- |x-y|, \phi_H^n(x,y)= \min(x+y - \frac{1}{n},2 - x- y + \frac{1}{n}), \phi_H^{\infty}(x,y)= \lim\limits_{n \rightarrow \infty} \phi_H^n(x,y)
\end{equation}
Finally, for any matrix $B$, let $\beta_h(B)$ denote the $h$th moment of its ESD.
Then the \emph{trace formula} implies
\begin{eqnarray}\label{trace_top}
\frac{1}{n}\Tr[BT_{n}]^h &=&
\frac{1}{n}\sum\limits_{1 \leq i_1,i_2,\ldots,i_h \leq n} \Bigg(\prod\limits_{j=1,2,\ldots,h-1} \frac{x_{L^{T}(i_j,
i_{j+1})}}{\sqrt{\phi_{T}(L^{T}(i_j,i_{j+1}))}}\Bigg) \times \frac{x_{L^{T}(i_h, i_1)}}{\sqrt{\phi_{T}(L^{T}(i_{h},i_{1}))}}\\
\E[\beta_{h}(BT_{n})] &=&
\E[\frac{1}{n}\Tr(BT_{n})^h] =
\frac{1}{n} \sum_{\pi: \ \pi \ \text{circuit}} \frac{\E
\X_{\pi}}{\prod\limits_{i=1,2,\ldots,h} \sqrt{\phi_T(\pi(i-1),\pi(i))}}.
\end{eqnarray}
\begin{eqnarray}\label{trace_hank}
\frac{1}{n}\Tr[BH_{n}]^h &=&
\frac{1}{n}\sum\limits_{1 \leq i_1,i_2,\ldots,i_h \leq n} \Bigg(\prod\limits_{j=1,2,\ldots,h-1} \frac{x_{L^{H}(i_j,
i_{j+1})}}{\sqrt{\phi_{H}(L^{H}(i_j,i_{j+1}))}}\Bigg) \times \frac{x_{L^{H}(i_h, i_1)}}{\sqrt{\phi_{H}(L^{H}(i_{h},i_{1}))}}\\
\E[\beta_{h}(BH_{n})] &=&
\E[\frac{1}{n}\Tr(BH_{n})^h] =
\frac{1}{n} \sum_{\pi: \ \pi \ \text{circuit}} \frac{\E
\X_{\pi}}{\prod\limits_{i=1,2,\ldots,h} \sqrt{\phi_{H}(\pi(i-1),\pi(i))}}.
\end{eqnarray}
\noindent
 {\bf Matched circuits:} Any value $L(\pi(i-1),\pi(i))$ is
an $L$ value of
$\pi$
 and $\pi$
 has an {\bf edge of order} $e\ (1 \leq e \leq h)$ if it has
  an \emph{L-value} repeated exactly $e$ times.
  If $\pi$ has at least one edge of
  order one then $\E(\X_{\pi})=0$. Thus only those $\pi$ with all $e \geq 2$ are relevant.
  Such circuits will be said to be
  \emph{matched}.
  $\pi$  is
  \textit{pair matched} if all its edges are of order two.

   \noindent
{\bf Equivalence relation on circuits:} Two circuits $\pi_1$ and $\pi_2$ are
equivalent iff their $L$-values agree at exactly the same pairs $(i,j)$. That is,
iff $\big\{L(\pi_1(i-1), \pi_1(i)) = L(\pi_1(j-1), \pi_1(j))$
$\Leftrightarrow L(\pi_2(i-1), \pi_2(i)) = L(\pi_2(j-1),
\pi(j))\big\}. $ This defines an equivalence relation between the
circuits.

\noindent {\bf Words:} Equivalence classes may be identified with partitions of $\{1,2,\cdots,h\}$:
to any partition we associate a \textbf{word} $w$ of length $l(w)=h$ of letters where  the first
occurrence of each letter  is in alphabetical order. For example, if $h = 6$, then the  partition
$\{ \{1,3,6\}, \{2,5\}, \{4\}\} $ is associated with $w=abacba$. For a word $w$, let $w[i]$
denote the letter in the $i$th position. The notion of matching and order $e$ edges
carries over to words. For instance, $abacabc$ is matched.
$abcadbaa$ is nonmatched, has
edges of order $1$, $2$ and $4$ and the corresponding partition is $\{\{1, 4, 7,
8\}, \{2, 6\}, \{3\}, \{5\}\}.$

\noindent
{\bf Independent vertex}: If $w[i]$ is the first occurrence of a letter then $\pi(i)$ is called a
\emph{independent vertex}. We make the convention that $\pi(0)$ is also an independent vertex. The other
vertices will be called \emph{dependent} vertices. If a word has $d$ distinct letters then there are
$d+1$ independent vertices.
\subsection{Reduction in the number of terms}\label{subsection:reductioninterms}
Fix an integer $h$. Define,
$$\Pi_{h}^{3+} = \{ \pi : \pi \ \text{is matched, of length}\  h \ \text{and has an edge of order greater than equal to } 3\}.
$$
$$S_{h}^{A_n} =\frac{1}{n}\sum_{\pi : \pi \in \Pi_{h}^{3+} } \frac{1}{ \prod\limits_{i=1}^h \sqrt{ \phi_A(\pi(i-1),\pi(i)) } }, \ \ A=H \ \ \text{or}\ \ T.$$
\begin{lem} \label{lem:Sh}$S_{h}^{A_n} \rightarrow 0$ as $n \rightarrow \infty$ for $A_n=T_n$ or $H_n$.  Hence, only pair matched circuits
are relevant while calculating $\lim E(\beta_h(A_n))$.
\end{lem}
\noindent \textbf{Proof.} We provide the proof only for $T_n$. Proof for $H_n$  is similar and details are omitted.
Note that,
$$S_{h}^{A_n} =\sum_w \frac{1}{n}\sum_{ \pi \in \Pi(w)\bigcap \Pi_{h}^{3+} } \frac{1}{ \prod_{i=1}^h \sqrt{ n - | \pi(i-1) - \pi(i) | }}=\sum_w S_{h,w} \ \ \text{say}.$$
It is enough to prove that for each $w$, $S_{h,w}\to 0$. We first restrict attention to
 $w$  which have only one edge of order 3 and all other edges of order 2.
 Note that this forces $h$ to be odd. Let $h
= 2t +1$ and $ |w| = t$.
Fix the $L$-values at
say $k_{1},k_{2},\ldots, k_{t}$
where $k_{1}$ is the
$L$-value corresponding to the order $3$ edge and let $i_0$ be such that $L(\pi(i_{0} -1) , \pi( i_{0} ) ) = k_{1} $. We start counting the number of possible $\pi$'s from the edge $(\pi(i_{0} -1) ,\pi(i_{0}) )$. Clearly the number of possible choices of
that edge is at most $2(n- k_{1})$. Having chosen the vertex $i_{0}$,
number of possible choices
of the vertex $(i_{0} + 1)$ is at most 2. Carrying on with this argument, we may  conclude that
the total number of $\pi$'s having $L$ values $k_{1},k_{2},\ldots,k_{t}$ is at most $C\times
(n -k_{1})$. Hence for some generic constant $C$,
\begin{eqnarray*} S_{h,w} = \frac{1}{n} \sum\limits_{k_i \in \{0,1,\ldots,n-1\} }\sum\limits_{{\pi: \pi \text{has $L$ values} }\atop {k_1,k_2,\ldots,k_t}} \frac{1}{(n-k_1)^{\frac{3}{2}} \prod\limits_{i=2}^t (n-k_i)}&\leq &
\frac{1}{n} \sum_{ k_{i} \in \{0,1,\ldots,n-1\} } \frac{ C \times ( n- k_{1}
) } { (n - k_{1}) ^{\frac{3}{2}} \prod\limits_{i=2}^t (n-k_i) } \\
&  = &\frac{O(\sqrt{n}) O ( (\log n)^{t-1} ) } {n} \rightarrow 0 \ \text{as}\ n \rightarrow \infty.
\end{eqnarray*}
In the last step, we have used the facts that \
 $ \sum_{k=1}^{n} \frac{1}{k} =O( \log{n})$  \  and for \ $0 < s <1$, $\sum_{k=1}^{n}
\frac{1}{k^{s}} = O(n^{1 -s })$.
It is easy to see that when $w$  contains more than one edge of order $3$ or more,  the order of the sum will be even smaller. This completes the proof of the first part.  The second part is immediate since $E(X_\pi)=1$ for every pair matched circuit and $E(|X_\pi|) < \infty$ uniformly over all $\pi$.
\qed
\subsection{Slope in balanced Toeplitz matrices}
\label{subsection: reduction_2}
Since the Toeplitz matrices have link function $L(i,j)=|i-j|$, given an $L$-value and a vertex there are at most two possible choices of the other vertex. Bryc, Dembo and Jiang  (2006)\cite{bryc}
showed that out of these two possible choices of vertices only one choice counts in the limit.
We show now that the same is true for the balanced matrices.
Let
$$\Pi_{h,+} = \{ \pi  \ \text{pair matched}: \ \text{there exists at least one pair} \ (i_0, j_0) \ \text{with} \
\pi(i_{0} -1) -\pi(i_{0}) + \pi(j_{0} -1) - \pi(j_{0}) \neq 0 \},$$
Define $\Pi_{h,+}(w)=\Pi_{h,+} \cap \Pi(w)$ and let $\pi(i-1) -\pi(i) $ be $i^{th}$ slope value,

\begin{lem}
\label{lem:opposite} Let $\pi \in \Pi_{h+} \; \text{and }  k_{1},k_{2},\ldots, k_{h}$ be the
$L$-values of $\pi$. Then $\exists \; j_{0} \in \{1,2,\ldots,h\} \; \text{such that } k_{j_{0}} = \Lambda(k_{1},k_{2},\ldots,k_{j_{0}
-1},k_{j_{0}+1},\ldots,k_{h} )$ for some linear function $\Lambda$.
\end{lem}
\noindent \textbf{Proof.} We note that sum of all the slope-values of $\pi$ is zero. Now sum of slope-values
from $j^{th}$ matched pair is  0 if the $L$ values are of opposite sign and $2k_{j}$ or $-2k_{j}$
if the $L$ values are of same sign. Hence we have, \ $f(k_{1},k_{2},\ldots,k_{h}) = 0$ for some linear function $f$ where coefficient of $k_{j} = 0$ if $L$ values corresponding to $k_{j}$ are of opposite sign and it is $\pm 2$ if $L$ values  corresponding to $k_{j}$ are of
same sign and the slope-values are positive (negative). Since $\pi \in \Pi_{h,+}, \ \exists \ k_{j} \neq 0$ such that $L$-values corresponding to $j^{th}$ pair have same sign. Let
$$\{ i_{1},i_{2},\ldots,i_{l} \} = \{ j : \; \text{coefficient of} \  k_{j} \neq 0 \} \ \ \text{and} \ \ j_{0} =\max \{ j : \ \text{coefficient of}\  k_{j} \neq 0 \}.$$
Then  $k_{j_{0}}$ can be expressed as
a linear combination
$k_{j_{0}} = \Lambda (k_{i_{1}},k_{i_{2}},\ldots,k_{i_{l}} ).$\qed
\begin{lem}\label{lem:Shplus}  \ $S_{h+} \stackrel{def}{=}
\frac{1}{n}\displaystyle{\sum_{\pi \in \Pi_{h,+}}} \frac{1}{ \prod_{i=1}^{h} \sqrt{ n- | \pi(i-1) -
\pi(i) | }} \rightarrow
0$ \ as $ n \rightarrow \infty$. Hence, to calculate $\lim E(\beta_h(BT_{n}))$ we may restrict attention to  pair matched circuits  where each edge has oppositely  signed $L$-value.
\end{lem}
\noindent
\textbf{Proof.} As in Lemma \ref{lem:Sh}, write $S_{h+} =\sum_w S_{h+, w}$
where $S_{h+, w}$ is the sum restricted to  $\pi \in \Pi_{h,+}(w)$.
Enough to show that this tends to zero for each $w$.
Let the corresponding $L$  values to this $w$ be $k_{1},k_{2},\ldots,
k_{h}$. Hence,
$$S_{h+,w} = \frac{1}{n}\displaystyle{\sum_{{k_{1},k_{2},\ldots,k_{h}} \atop{ \in \{ 0,1,2,\ldots,n-1\} }}}
 \frac{(\# \ \text{of}\
 \pi \in \Pi_{h+}(w) \; \ni \; L \ \text{values of}\ \pi \ \text{are}\  \{ k_{1},k_{2},\ldots,k_{h} \} )} { \displaystyle{\prod_{i=1}^h}( n- k_{i} )  }.$$
For this  fixed set of
$L$ values,
 there are at most $2^{2h}$ sets of slope-values.
 It is enough to prove the result for any one such set.
  Now we start counting the number of possible $\pi$'s having those
slope values.

\noindent
 By the previous lemma there exists $j_{0}$ such that $k_{j_{0}} = \Lambda (
k_{i_1},k_{i_2},\ldots,k_{i_l} )$. We start counting number of possible $\pi$ from the edge
corresponding to the $L$ value $k_{j_0}$, say, $( \pi(i_{*} -1) ,\pi(i_{*}) )$. Clearly
number of ways to choose vertices $\pi(i_{*} -1)$ and $\pi(i_{*} ) $ is $ (n- k_{j_0})$. Having
chosen $\pi(i_{*})$,
there is only one choice of $\pi(i_{*} +1)$ (since the slope-values have been fixed),
We continue this procedure to choose all the vertices of the circuit $\pi$ and hence number of $\pi$'s having the fixed set of slope-values is at most $(n-k_{j_0} )$. Note that since $w$ and the slope signs are fixed, the linear function $\Lambda$ and the index $j_0$ are determined as well. Thus for that fixed set we have,
\begin{eqnarray*}
S_{h+,w}^{set}  &\leq &\frac{1}{n}
\sum\limits_{k_i  \in \{0,1,2,\ldots,n-1 \}}
\frac{n - k_{j_0} } { \prod_{i=1}^h(n-k_i) }
 =  \frac{1}{n} \sum\limits_{ {k_i \in \{ 0,1,2,\ldots ,n-1\}} \atop { i\neq j_0}}
\frac{1}{\prod_{{i=1} \atop {i \neq j_0}}^h (n-k_i)}
\end{eqnarray*}
As $k_{j_0} = \Lambda(k_{i_1},k_{i_2},\ldots,k_{i_l} )$, in
the above sum, $k_{j_0}$ should be kept fixed which implies that $ S_{h+,w} \leq
\frac{O\left((\log n)^{h-1} \right) } {n} \rightarrow 0 \ \text{as} \  n \rightarrow \infty$, proving the first part. The second part now follows immediately.
\qed
\subsection{Convergence of the moments $\E[\beta_h(BT_n)]$ and $\E[\beta_h(BH_n)]$}
\label{subsection:conv}
We need first establish a  few results on  moments of truncated uniform.
For a given random variable $X$ (to be chosen), define (whenever it is finite)
$$g_T(x) = \E [ \phi_T^n(X,x) ^{-(1+\alpha)} ] \text{ and } g_H(x)= \E[\phi_H^n(X,x)^{-(1+ \alpha)}],$$
\begin{lem}
\label{lem: g(x)} Let $ x \in \N_{n}= \{1/n,2/n,\ldots,1\},\ \alpha > 0 \ \text{and}  \ X \text{ be
discrete uniform on} \ \ \ \N_{n}$. Then, for some  constants $C_{1},C_{2}$,
$$\max\{g_T(x),g_H(x) \} \leq C_{1} x^{-\alpha} + C_{2} (1
-x+1/n) ^{-\alpha} +1/n. $$
\end{lem}
\noindent \textbf{Proof.} Note that
$$g(x)=\frac{1}{n}\sum\limits_{y=1}^{n} \frac{1}{[1 - | x- \frac{y}{n} |]^{1+\alpha}}
=\frac{1}{n}\sum\limits_{ y<j} \frac{1}{(1 - \frac{j -y}{n} )^{1+\alpha} } + \frac{1}{n}\sum\limits_{y >j}
\frac{1}{(1- \frac{y-j}{n})^{1+\alpha}} + \frac{1}{n},$$
where $x=j/n$ and $1 < j <n$. For $j=1$ or $n$ similar arguments will go through.
Now,
\begin{eqnarray*}
\frac{1}{n}\sum\limits_{y<j} (1 -\frac{j-y}{n} ) ^{ -(1+\alpha)}  = \frac{1}{n} \sum\limits_{t=1}^{j-1} (1 - \frac{t}{n} )^{-(1+\alpha)}
& =&n^{\alpha} \sum\limits_{t=n-j+1}^{n-1} t^{-(1+\alpha)}  \\
&\leq &  n^{\alpha} \times \frac{C_1}{(n-j+1)^{\alpha}}=   C_{1} (1 - x+ 1/n)^{-\alpha}.
\end{eqnarray*}
By similar arguments,
$\frac{1}{n} \sum\limits_{y >j} \frac{1}{(1 - \frac{y-j}{n} ) ^{1+\alpha}} \leq C_{2} x^{-\alpha}.$

\noindent
By similar calculations $g_H(x)  \leq C_{1} x^{-\alpha} + C_{2} (1
-x+1/n) ^{-\alpha} +1/n $ and thus the result follows.
\qed
\begin{lem}
\label{lem: x^a} Suppose $U_{i, n}$ are i.i.d. discrete uniform on  $\N_{n}.$
Let $a_{i} \in \Z, 1 \leq i \leq m$ be fixed and $0<\beta<1.$ Let, $Y_n =  \sum\limits_{i=1}^m a_i
U_{i,n}$ and $Z_n = 1- Y_n +1/n$. Then

\[ \sup_{n} \E |Y_n| ^{-\beta} I (|Y_n| \geq 1/n)\  +\
 \sup_{n} \E |Z_n| ^{-\beta} I (|Z_n| \geq 1/n) < \infty.\]
\end{lem}
\noindent \textbf{Proof.}  First note that
\begin{eqnarray*}
\pr( |Y_{n}| \leq M/n ) & = & \E \Big[\pr\Big( -M/n \leq \sum\limits_{i=1}^{m} a_{i} U_{i,n} \leq M/n \Big| U_{i,n}, j\neq i_{0} \Big) \Big] \\
& = & \E \Big[ \pr \Big( -M/n -\sum\limits_{{i=1} \atop {i \neq i_{0}} } ^{m} a_{i} U_{i,n} \leq a_{i_{0}} U_{i_{0},n} \leq M/n -\sum\limits_{{i=1} \atop {i \neq i_{0}}}^{m} a_{i} U_{i,n}\Big) \Big] \leq  (2M+1)/n.
\end{eqnarray*}
Let $U_1,U_2,\ldots,U_m$ be $m$ i.i.d $U(0,1)$ random variables. We note that,
$$\Big( U_{1,n}
,U_{2,n} ,\ldots,U_{m,n} \Big) \stackrel{\mathcal{D}}{=}   \Big( \frac{ \lceil n U_1\rceil}{n} ,
\frac{\lceil n U_2 \rceil} {n} ,\ldots, \frac{\lceil n U_m \rceil}{n}\Big).$$ Define
$$\hat{Y}_n =\sum\limits_{i=1}^{m} a_{i} \frac{ \lceil n U_i \rceil}{n} \ \ \text{,} \ \ Y =\sum\limits_{i=1}^{m} a_i U_i
\ \text{and}\  K=\sum\limits_{i=1}^{m} |a_{i}|.$$ Then,
$$\hat{Y}_n \stackrel{\mathcal{D}}{=} Y_n \ \ \text{and} \ \  |\hat{Y}_n - Y | \leq K/n.$$
\begin{eqnarray*}
 \E |Y_n| ^{-\beta} I (|Y_n| \geq 1/n) & = & \E |Y_n| ^{-\beta} I (1/n \leq |Y_n| \leq 2K/n)  +   \E |\hat {Y}_n| ^{-\beta} I (| \hat {Y}_n| > 2K/n) \\
 & \leq & n^{\beta} \frac{4K+1}{n} +  \E (|Y| - K/n)  ^{-\beta} I (| \hat {Y}_n| > 2K/n) \\
  & \leq & o(1)+  \E (|Y| - K/n)  ^{-\beta} I (| Y| > K/n)\\
    & \leq & o(1)+  \int_{ x > K/n} (x - K/n)^{-\beta} f_{|Y|}(x) dx.
\end{eqnarray*}
Now,
$$\int_{ x > K/n} (x - K/n)^{-\beta} f_{|Y|}(x) dx = \int_{ 0}^{\infty} x ^{-\beta} f_{|Y| }(x+ K/n) dx.$$
It is easy to see that $f_{Y}$ vanishes outside $[-K,K]$. Using induction one can also prove that,
$f_Y(x) \leq 1\  \forall \ \ x$. These two facts yields,
\[ \int_{ 0}^{\infty} x ^{-\beta} f_{|Y| }(x+ K/n) dx \le \int_{ 0}^{K+ K/n} x ^{-\beta} 2 dx  =O(1). \]
Hence,
$$\sup\limits_{n} \E |Y_n|^{-\beta} I ( |Y_n| \geq 1/n ) < \infty. $$
The proof of the finiteness of the other supremum is similar and we omit the details.  \qed
\begin{lem}
\label{lem:balancedm1} Suppose $\{ x_i \}$ are i.i.d. bounded with mean zero and variance 1.
Then $\lim \E[\beta_h(BT_n)]$ and $\E[\beta_h(BH_n)]$ exists for every $h$.
\end{lem}
\noindent
\textbf{Proof.} From Lemma \ref{lem:Sh} it follows that $\E[\beta_{2k+1}(BA_n)]\rightarrow 0$ as $n \rightarrow \infty$ where $A_n=T_n$ or $H_n$. From Lemma \ref{lem:Sh} and Lemma \ref{lem:Shplus} (if limit exists) we have,
\begin{eqnarray*}
\lim\limits_{n \rightarrow \infty} \E[\beta_{2k}(BT_n)] &= & \sum\limits_{w \text{ pair matched}} \lim\limits_{n \rightarrow \infty} \frac{1}{n} \sum\limits_{ \pi \in \Pi^{*}(w)} \frac{\E \X_{\pi}}{\prod_{i=1}^h \sqrt{\phi_T(\pi(i-1),\pi(i))}}\\
& = & \sum\limits_{w \text{ pair matched}} \lim\limits_{n \rightarrow \infty} \frac{1}{n} \sum\limits_{ \pi \in \Pi^{**}(w)} \frac{1}{\prod_{i=1}^h \sqrt{\phi_T(\pi(i-1),\pi(i))}}
\end{eqnarray*}
where $\Pi^{**}(w)= \{ \pi: \ w[i]=w[j] \Rightarrow \pi(i-1)-\pi(i)+ \pi(j-1) - \pi(j)=0 \}$. Denote $x_i = \pi(i)/n$. Let $S= \{0\} \cup \{ \min(i,j): w[i]=w[j], i \neq j \}$ be the set of all independent vertices of the word $w$ and let $x_S= \{x_i: i \in S\}$. Each $x_i$ can be expressed as a unique linear combination $L_i^T(x_S)$.  $L_i^T$ depends on word $w$ but for notational convenience we suppress its dependence. Note that $L_i^T(x_S)=x_i$ for $i \in S$ and also summing $k$ equations we get $L_{2k}^T(x_S)=x_0$. If $w[i]=w[j]$ then $|L_{i-1}^T(x_S) - L_i^T(x_S)| = |L_{j-1}^T(x_S) - L_j^T(x_S)|$. Thus using this equality and proceeding as in Bose and Sen \cite{bosesen08} and [BDJ] \cite{bryc} we have,
$$\lim\limits_{n \rightarrow \infty} \E[\beta_{2k}(BT_n)]=  \sum\limits_{w \text{ pair matched}} \lim\limits_{n \rightarrow \infty} \E\bigg[\frac{\I(L_i^T(U_{n,S}) \in \N_n, i \notin S \cup \{2k\})}{\prod\limits_{i \in S\setminus \{0\}} \phi_T^n(L_{i-1}(U_{n,S}),U_i))}\bigg],$$
where for each $i \in S$, $U_{n,i}$ is discrete uniform on $\N_n$ and $U_{n,S}$ is the random vector on $\R^{k+1}$ whose co-ordinates are $U_{n,i}$ and $U_{n,i}$'s are independent of each other.
We claim that
\begin{equation}\label{eq:limexptoep}
\lim\limits_{n \rightarrow \infty} \E[\beta_{2k}(BT_n)]=  m_{2k}^T=\sum\limits_{w \text{ pair matched}}  m_{2k,w}^T= \sum\limits_{w \text{ pair matched}}  \E\bigg[\frac{\I(L_i^T(U_{S}) \in (0,1), i \notin S \cup \{2k\})}{\prod\limits_{i \in S\setminus \{0\}} \phi_T^{\infty}(L_{i-1}(U_{S}),U_i))}\bigg],
\end{equation}
where for each $i \in S$, $U_i \sim U(0,1)$ and $U_S$ is a $\R^{k+1}$ dimensional random vector whose co-ordinates are $U_i$ and they are independent of each other. Note that to prove (\ref{eq:limexp}) it is enough to show that for each pair matched word $w$ and for each $k$ there exists $\alpha_k>0$ such that
\begin{equation}
\sup\limits_n \E\bigg[\bigg(\frac{\I(L_i^T(U_{n,S}) \in \N_n, i \notin S \cup \{2k\})}{\prod\limits_{i \in S\setminus \{0\}} \phi_T^n(L_{i-1}(U_{n,S}),U_i))}\bigg)^{1+\alpha_k}\bigg] < \infty
\end{equation}
We will prove that for each pair matched word $w$
$$\sup\limits_n \E\bigg[\bigg(\frac{\I(L_i^T(U_{n,S}) \in \N_n, i \notin S \cup \{2k\}, i < \max S)}{\prod\limits_{i \in S\setminus \{0\}} \phi_T^n(L_{i-1}(U_{n,S}),U_{n,i}))}\bigg)^{1+\alpha_k}\bigg] < \infty $$  and we prove it by induction on $k$. For $k=1$ the expression reduces to \ $\E\Big[ \Big( \frac{1}{1 -|U_{n,0} - U_{n,1}|} \Big) ^{1+\alpha} \Big].$
Now,
\begin{eqnarray*}
\E\Big[ \Big( \frac{1}{1 -|U_{n,0} - U_{n,1}|} \Big) ^{1+\alpha} \Big] & = & \E\Big[ \E\Big\{ \Big( \frac{1}{1 -|U_{n,0} - U_{n,1}|} \Big) ^{1+\alpha} \Big| U_{n,0} \Big\}\Big]\\
& = & \E [ g_T(U_{n,0})] \leq  C_1\E[U_{n,0}^{-\alpha}] + C_2\E[(1 - U_{n,0})^{-\alpha}]+ 1/n \ \ \ \text{ by Lemma
\ref{lem: g(x)} }.
\end{eqnarray*}
Now by Lemma \ref{lem: x^a} we have, $\sup\limits_{n} \E\Big[ \Big( \frac{1}{1 - |U_{n,0} - U_{n,1}|} \Big)^{1+\alpha} \Big] < \infty \ \ \forall \ \ 0<\alpha<1.$

\noindent
Now we assume the result for $k=1,2,\ldots,t.$ We prove it for $k=t+1$. Fix any pair matched word $w_0$. Note that the random variable corresponding to the generating vertex of the last letter appears only once and hence we can do the following calculations. Let
$$B_{t+1}= \Bigg[\frac{I\Big(L_{i}^{T}(U_{n,S}) \in \N_{n}, i \notin S \cup \{2(t+1)\}, i < \max S\Big)}{ \prod\limits_{i \in S \setminus \{0\}}(1- |L_{i-1}^T(U_{n,S})-U_{n,i}|)}\Bigg]^{1+\alpha}. $$ Then
\begin{eqnarray*}
\E [B_{t+1}] & = &\E\Big[ \E[B_{t+1} \Big| U_{n,i}, i \in S\setminus\{i_{t+1}\}] \Big]\\
& = &\E\Bigg[ \underbrace{\Bigg(\frac{I\Big(L_{i}^{T}(U_{n,S}) \in \N_{n}, i \notin S , i < \max S\setminus\{i_{t+1}\}\Big)}{ \prod\limits_{i \in S \setminus \{0,i_{t+1}\}}(1- |L_{i-1}^T(U_{n,S})-U_{n,i}|)}\Bigg)^{1+\alpha} }_{\Phi_n}\times \underbrace{g_T(U_{n,i_{t+1} -1}) I [L_{i_{t+1}-1}^{T}(U_{n,S}) \in \N_{n}] }_{\Psi_n}\Bigg]
\end{eqnarray*}
By Lemma \ref{lem: g(x)} and Lemma \ref{lem: x^a}, $\sup\limits_n ||\Psi_n||_q < \infty$ whenever $\alpha q <1$. Here $||\cdot||_q$ denotes the $L_q$ norm.

Let us now consider the word $w_{0}^{*}$ obtained from $w_{0}$ by removing both
occurrences of the last used letter. We note that the quantity
$\Phi_{n}$ is the candidate for the
expectation expression corresponding to the word $w_{0}^{*}$. Now by induction hypothesis, $\exists \ \ \alpha_{t} > 0$ such that,
$$\sup\limits_{n}\E\bigg[\Bigg(\frac{ I \Big(L_{i}^{T}(U_{n,S}) \in \N_{n}, i \notin S i < \max S\setminus \{i_{t+1}\}\Big) }{\prod\limits_{i \in S \setminus \{0,i_{t+1}\}} (1- |L_{i_1}^T(U_{n,S})-U_{n,i}|)} \Bigg)^{1+\alpha_{t}}\bigg] < \infty.$$
Hence $\sup\limits_{n} ||\Phi_n||_p < \infty \ \ \text{if, } (1+\alpha)p \leq (1+\alpha_{t}).$
Therefore
$$\alpha_{t+1} + \frac{1+\alpha_{t+1}}{1+\alpha_{t}} < \frac{1}{p} + \frac{1}{q} = 1 \Rightarrow
 \sup\limits_n \E[\Phi_n \Psi_n] \leq \sup\limits_n||\Phi_n||_p ||\Psi_n||_q < \infty.$$
This proves the claim for balanced Toeplitz matrix. For balanced Hankel matrix we again use Lemma  \ref{lem: g(x)} and Lemma \ref{lem: x^a} and proceed in an exact similar way to get,
\begin{equation}\label{eq:limexphan}
\lim\limits_{n \rightarrow \infty} \E[\beta_{2k}(BH_n)]=  m_{2k}^H=\sum\limits_{{w \text{ pair matched}} \atop {\text{and symmetric}} } m_{2k,w}^H=\sum\limits_{{w \text{ pair matched}} \atop {\text{and symmetric}} } \E\bigg[\frac{\I(L_i^H(U_{S}) \in (0,1), i \notin S \cup \{2k\}}{\prod\limits_{i \in S\setminus \{0\}} \phi_H^{\infty}(L_{i-1}(U_{S}),U_i))}\bigg].
\end{equation}
Symmetric pair matched words are those in which every letter appears once each in an odd position and  an even position. Using ideas of Bose and Sen (2008) \cite{bosesen08} and [BDJ] \cite{bryc} it can be shown that for any pair matched non-symmetric word $w$, $m_{2k,w}^H=0$. So the above summation is taken over only pair matched symmetric words.
\qed
\subsection{An approximation result}\label{subsection:approximation}
Even though the limit of the moments have been established, it does not seem to be easy to show that this moment sequence determines a probability distribution uniquely (which would then be the candidate  LSD). We tackle this issue by using approximating matrices whose scalings are not unbounded.
We shall use the
L\'evy distance metric to develop this approximation. Recall that this metric metrizes weak convergence of probability measures on ${\mathbb R}$.
Let $\mu_i$, $i=1, 2$be two probability measures on ${\mathbb R}$. The L\' evy distance between them
is given by,
$$ \rho(\mu_1,\mu_2)= \inf \{ \varepsilon > 0: \ F_1(x - \varepsilon) - \varepsilon < F_2 (x) < F_1(x +\varepsilon) + \varepsilon, \ \forall \ x \in {\mathbb R} \},$$
where $F_i$ $ i=1, 2$ are the distribution functions corresponding to the measures $\mu_i, i =1, 2$.
\begin{prop} \label{prop:trunc}
(Bhamidi, Evans and Sen (2009))
Suppose $A_{n\times n}$ is  real symmetric matrix and $B_{m \times m}$ is the principal sub-matrix of $A_{n\times n}$. Then
$$\rho(F^A,F^B) \leq \Big(\frac{n}{m}-1\Big) \wedge 1.$$
Let $(A_{k})_{k=1}^{\infty}$ be a sequence of real symmetric matrices.
For each $\varepsilon > 0$, and each $k$, let $(B_K^{\varepsilon})_{k=1}^{\infty}$ be
an $n_k^{\varepsilon} \times n_k^{\varepsilon}$ principal sub-matrix of $A_k$.
Suppose that for each $\varepsilon > 0$, $F^{\varepsilon}_{\infty}= \lim\limits_{k \rightarrow \infty} F^{B_k^{\varepsilon}}$ exists and $\limsup\limits_{k \rightarrow \infty} n_k/n_k^{\varepsilon} \leq 1+ \varepsilon$. Then $F_{\infty} = \lim\limits_{k \rightarrow \infty} F^{A_k}$ exists and is given by $F_{\infty}= \lim\limits_{\varepsilon \downarrow 0} F^{\varepsilon}_{\infty}$.
\end{prop}
\noindent
    Consider the principal submatrix $BT_{n}^{\varepsilon}$ of $BT_n$ obtained by retaining the first  $n(1-\varepsilon)$ rows and columns of $BT_n$.
    Then for this matrix, since  $|i-j| \leq n(1 - \varepsilon)$, the balancing factor becomes bounded.
We shall show that LSD of  $\{F^{BT_n^{\veps}}\}$ exists  for every $\varepsilon$ and then invoke the above
result to obtain the LSD  of $\{F^{BT_n}\}$. Similar argument holds for $\{BH_n\}$, by considering the principal sub-matrix obtained by removing the first $n \veps/2$ and last $n \veps/2$ rows and columns.
\subsection{Existence of limit of $\{F^{BT_n^{\veps}}\}$ and $\{F^{BH_n^{\veps}}\}$ almost surely}\label{subsection:limitoftruncation}
Clearly, for any fixed $\varepsilon >0$, we may write,
\begin{equation}
\frac{1}{n}\Tr[BT_{n}^{\varepsilon}]^h= \frac{1}{n}\sum\limits_{\pi: \pi \text{ circuit}}\frac{\X_{\pi}}{\prod\limits_{i=1,2,\ldots,h} \sqrt{n - |\pi(i-1)-\pi(i)|}}\times \prod\limits_{i=1,2,\ldots,h} \I\Big[\pi(i) \leq n(1- \varepsilon)\Big]
\end{equation}
\noindent
And similarly,
\begin{equation}
\frac{1}{n}\Tr[BH_{n}^{\varepsilon}]^h= \frac{1}{n}\sum\limits_{\pi: \pi \text{ circuit}}\frac{\X_{\pi}}{\prod\limits_{i=1,2,\ldots,h} \phi_H(\pi(i-1),\pi(i))}\times \prod\limits_{i=1,2,\ldots,h} \I\Big[n \veps/2 \leq \pi(i) \leq n(1- \varepsilon/2)\Big]
\end{equation}
Since for every $\veps >0$, the scaling is bounded, the proof of the following Lemma is exactly as the proof of Lemma 1 and Theorem 6 of Bose and Sen (2008) \cite{bosesen08}. Hence we skip the proof.
\begin{lem}\label{lem:epsm1} (i) If $h$ is odd, $\E[\beta_{h}(BT_n^{\veps})]  \rightarrow 0$ and $\E[\beta_{h}(BH_n^{\veps})]  \rightarrow 0$.

 (ii) If $h$ is even ($=2k$), then (below the sums are over all pair matched $w$)
  $$\lim\limits_{n \rightarrow \infty} \E[\beta_h(BT_n^{\veps})]= \sum_{w}
   p_{BT^{\veps}}(w)
       =  \underbrace{\int_0^{1-\veps} \cdots \int_0^{1-\veps}}_{k+1}  \frac{\prod\limits_{i \notin S \cup \{2k\}} \I (0 \leq L_i^T(x_S) \leq 1- \veps)}{ \prod\limits_{i \in S\setminus \{0\}} (1 - |L_{i-1}^T(x_S) - x_i|)} dx_S =\  \beta_{2k}^{T^\veps} \ \ \text{say}
        $$
            $$\lim\limits_{n \rightarrow \infty} \E[\beta_h(BH_n^{\veps})]= \sum_{w}
      p_{BH^{\veps}}(w)=
           \underbrace{\int_{\veps/2}^{1-\veps/2}  \cdots \int_{\veps/2}^{1-\veps/2}}_{k+1}  \frac{\prod\limits_{i \notin S \cup \{2k\}} \I (\veps/2 \leq L_i^H(x_S) \leq 1- \veps/2)}{ \prod\limits_{i \in S\setminus \{0\}} \phi_H^{\infty}(L_{i-1}^H(x_S),x_i)} dx_S= \beta_{2k}^{H^\veps} \ \ \text{say}.$$
       Further, $\max\{\beta_{2k}^{T^\veps}, \beta_{2k}^{H^\veps}\} \leq \frac{2k!}{k! 2^k} \times \veps^{-k}$. Hence
        there exists unique probability distributions $F^{T^\veps}$  and $F^{H^\veps}$ with $\beta_{k}^{H^\veps}$  and $\beta_{k}^{H^\veps}$ (respectively) as their moments.
\end{lem}
The almost sure convergence of $\{F^{BT_n^\veps}\}$ and $\{F^{BH_n^\veps}\}$ now follows from the following Lemma. We omit its proof since it is essentially a repetition of arguments Proposition 4.3 and Proposition 4.9  of[BDJ] who established it for the usual Toeplitz matrix $T_n$ and usual Hankel matrix $H_n$.
\begin{lem}\label{lem:epsalmost}
Fix any $\veps>0$ and let $A_n=T_n$ or $H_n$. If the input sequence is uniformly bounded, independent, with mean zero and variance one then
\begin{equation}
 \E\Big[ \frac{1}{n} \Tr(BA_n^{\veps})^h - \E \frac{1}{n} \Tr(BA_n^{\veps})^h\Big]^4 = O\Big(\frac{1}{n^2}\Big).
\end{equation}
As a consequence, the ESD of $BA_n^\veps$ converges to $F^{A^\veps}$ almost surely.
\end{lem}
\subsection{Connecting limits of $BT_n^\veps$ and $BT_n$}\label{subsection:connect}
From Lemma \ref{lem:epsm1} and Lemma \ref{lem:epsalmost}, given  any $\veps >0$, there exists $B_{\veps}$ such that $\pr(B_{\veps})=1$ and on $B_{\veps}$,  $F^{BT_n^{\veps}}\Rightarrow F^{T^\veps}$.

\noindent
Fix any sequence $\{\veps_m\}_{m=1}^{\infty}$ decreasing to $0$.
Define $B= \cap B_{\veps_m}$. Using Proposition \ref{prop:trunc}, on $B$, $F^{BT_n} \Rightarrow F^T$ for some non-random distribution function $F^T$ where, $F^T$ is the weak limit of $\{F^{T^{\veps_m}}\}_{m=1}^{\infty}$.

Let  $X^{\veps_m}$ (resp. $X$)  be a random variable with distribution  $F^{T^{\veps_m}}$ (resp. $F^T$) with $k$th moments  $\beta_{k}^{T^{\veps_m}}$  (resp. $\beta_k^T$).
From Lemma \ref{lem:epsm1},  and (\ref{eq:limexptoep})
it is clear that for all $ k \geq 1$,
$$\lim\limits_{m \rightarrow \infty} \beta_{2k+1}^{T^{\veps_m}}=0 \ \ \text{and} \ \ \lim\limits_{m \rightarrow \infty} \beta_{2k}^{T^{\veps_m}} =m_{2k}^T=\sum\limits_{w \text{ pair matched}}m_{2k,w}^T.$$
From Lemma \ref{lem:balancedm1}, $m_{2k}^T$ is finite for every $k$.
Hence $\{(X^{\veps_m})^k\}_{m=1}^{\infty}$ is uniformly integrable for every $k$ and  $\lim\limits_{m \rightarrow \infty} \beta_{k}^{T^{\veps_m}}= \beta_k^T$. This proves that $m_k^T=\beta_k^T $ and so $\{m_k^T\}$ are the moments of $F^T$.
The argument for $BH_n$ is exactly same and hence details are omitted. The proof of Theorem \ref{thm:main theorem} is now complete. \qed

\noindent \textbf{Acknowledgement}. We learnt the idea of the truncation that we have used, from personal communication with Arnab Sen who credits them to Steven Evans.

\end{document}